\documentclass[reqno,12pt]{amsart}

\usepackage{amssymb}
\newtheorem{Proposition}{Proposition}[section]
\newtheorem{Definition}[Proposition]{Definition}
\newtheorem{Lemma}[Proposition]{Lemma}
\newtheorem{Theorem}[Proposition]{Theorem}

\newtheorem{Corollary}[Proposition]{Corollary}
\newcommand{\R}{\mathbb{R}}
\DeclareMathOperator{\Val}{Val}
\DeclareMathOperator{\cnc}{cnc}
\DeclareMathOperator{\nc}{nc}
\DeclareMathOperator{\spt}{spt}
\DeclareMathOperator{\vol}{vol}
\DeclareMathOperator{\Gr}{Gr}
\DeclareMathOperator{\kl}{Kl}
\newcommand\Hom{\mathbf{Hom}}
\newcommand\PD{\mathbf{PD}}

\begin{document}

\title[Valuations and quaternionic integral geometry]{A product formula for valuations on manifolds with applications to the integral geometry of the quaternionic line}
\author{Andreas Bernig}
\email{andreas.bernig@unifr.ch}

\address{D\'epartement de Math\'ematiques, Chemin du Mus\'ee 23, 1700 Fribourg, Switzerland}
\begin{abstract}
The Alesker-Poincar\'e pairing for smooth valuations on manifolds is
expressed in terms of the Rumin differential operator acting on the
cosphere-bundle. It is shown that the derivation operator, the
signature operator and the Laplace operator acting on smooth
valuations are formally self-adjoint with respect to this pairing. As
an application, the product structure of the space of $SU(2)$- and
translation invariant valuations on the quaternionic line is
described. The principal kinematic formula on the quaternionic line $\mathbb{H}$ is stated and proved. 
\end{abstract}

\thanks{ 
 \\ Supported
  by the Schweizerischer Nationalfonds grants PP002-114715/1 and SNF 200020-105010/1.}
\maketitle 
\section{Smooth valuations on manifolds}

Let $M$ be a smooth manifold of dimension $n$. For simplicity, we suppose that $M$ is oriented, although the whole theory works in the non-oriented case as well. 
Following Alesker, we set $\mathcal{P}(M)$ to be the set of compact submanifolds with corners. 
\begin{Definition}
A valuation on $M$ is a real valued map $\mu$ on $\mathcal{P}(M)$ which is additive in the following sense:
whenever $X,Y,X\cap Y$ and $X \cup Y$ belong to $\mathcal{P}(M)$, then 
\begin{displaymath}
\mu(X \cup Y)+\mu(X \cap Y)=\mu(X)+\mu(Y). 
\end{displaymath}
\end{Definition}

A set $X \in \mathcal{P}(M)$ admits a conormal cycle $\cnc(X)$, which
is a compactly supported Legendrian cycle on the cosphere bundle
$S^*M$. Sometimes it will be convenient to think of $S^*M$ as the set of pairs
$(p,P)$ with $p \in M$ and $P \subset T_pM$ an oriented hyperplane, at
other places it is better to think of it as the set of pairs
$(p,[\xi])$ where $p \in M$ and $\xi \in T_p^*M\setminus \{0\}$ and the brackets
denote the equivalence class for the relation $\xi_1 \sim \xi_2 \iff
\xi_1=\lambda \xi_2, \lambda>0$. 

A valuation $\mu$ on $M$ is called {\it smooth} if there exist an $n-1$-form $\omega \in \Omega^{n-1}(S^*M)$ and an $n$-form $\phi \in \Omega^n(M)$ such that 
\begin{equation} \label{smooth_vals}
\mu(X)=\cnc(X)(\omega)+\int_X \phi, \quad X \in \mathcal{P}(X).
\end{equation}

If $\mu$ can be expressed in the form \eqref{smooth_vals}, we say that
$\mu$ {\it is represented} by $(\omega,\phi)$. The space of smooth valuations on $M$ is denoted by $\mathcal{V}^{\infty}(M)$. It is a Fr\'echet space (see \cite{ale05b}, Section 3.2 for the definition of the topology). If $M=V$ is a vector space, the subspace of translation invariant smooth valuations will be denoted by $\Val^{sm}(V)$. 

Let $N$ be another oriented $n$-dimensional smooth manifold and $\rho:N \to M$ an orientation preserving immersion. 
Then $\rho$ induces a map $\tilde \rho:S^*N \to S^*M$, sending $(p,P)$ to $(\rho(p),T_p\rho(P))$. It clearly satisfies $\pi \circ \tilde \rho=\rho \circ \pi$. 

The valuation $\rho^*\mu$ on $N$ such that 
\begin{displaymath}
\rho^* \mu(X):=\mu(\rho(X)), \quad X \in \mathcal{P}(N)
\end{displaymath}
is again smooth. If $\mu$ is represented by $(\omega,\phi)$, then $\rho^*\mu$ is represented by $(\tilde \rho^*\omega,\rho^*\phi)$. This follows from the fact that $\cnc(\rho(X))=\tilde \rho_* \cnc(X)$. Note also that $\widetilde{\rho^{-1}}=(\tilde \rho)^{-1}$ if $\rho$ is a diffeomorphism. 

We will use some results of \cite{bebr07}, which we would like to recall. 
The cosphere bundle $S^*M$ is a contact manifold of dimension $2n-1$
with a global contact form $\alpha$ ($\alpha$ is not unique,
but this will play no role here). The projection from $S^*M$ to $M$
will be denoted by $\pi$, it induces a linear map $\pi_*$ (fiber integration)
on the level of forms. 

Given an $n-1$-form $\omega$ on $S^*M$, there exists a unique vertical form $\alpha \wedge \xi$ such that $d(\omega+\alpha \wedge \xi)$ is vertical (i.e. a multiple of $\alpha$). The Rumin differential operator $D$ is defined as $D\omega:=d(\omega+\alpha \wedge \xi)$ \cite{rum94}. The following theorem was proved in \cite{bebr07}.

\begin{Theorem} \label{bernig_broecker}
Let $\omega \in \Omega^{n-1}(S^*M)$, $\phi \in \Omega^n(M)$ and define the smooth valuation $\mu$ by \eqref{smooth_vals}. Then $\mu=0$ if and only if 
\begin{enumerate}
\item $D\omega+\pi^*\phi=0$ and
\item $\pi_* \omega=0$ for all $p \in M$. 
\end{enumerate}
Moreover, if $D\omega+\pi^*\phi=0$, then $\mu$ is a multiple of the
Euler characteristic $\chi$. 
\end{Theorem}

The support of a smooth valuation $\mu$ is defined as 
\begin{displaymath}
\spt \mu:=M \setminus \left\{p \in M: \exists p \in U \subset M \text{
    open }, \mu|_U=0 \right\}. 
\end{displaymath}
The subspace of $\mathcal{V}^\infty(M)$ consisting of compactly
supported valuations will be denoted by $\mathcal{V}_c^\infty(M)$.

Let $\int: \mathcal{V}^\infty_c(M) \to \R$ denote the integration
functional \cite{ale05d}. If $\mu$ has compact support, then $\int
\mu:=\mu(X)$, where $X \in \mathcal{P}(M)$ is an $n$-dimensional
manifold with boundary containing $\spt \mu$ in its interior. It is
clear that, if $\mu$ is represented by $(\omega,\phi)$ with compact
supports, then 
\begin{displaymath}
\int \mu= \int_M \phi= [\phi] \in H^n_c(M)=\R. 
\end{displaymath} 

Before stating our main theorem we have to recall two other
constructions of Alesker. 

The first one is the Euler-Verdier involution $\sigma:
\mathcal{V}^\infty(M) \to \mathcal{V}^\infty(M)$ \cite{ale05b}. Let $s:S^*M \to S^*M$ be the natural
involution on $S^*M$, sending $(p,P)$ to $(p,\bar P)$, where $\bar P$
is the hyperplane $P$ with the reversed orientation. If a valuation $\mu \in \mathcal{V}^\infty(M)$ is represented by the
pair $(\omega,\phi)$, then $\sigma \mu$ is defined as the valuation
which is represented by the pair $((-1)^n s^*\omega,(-1)^n \phi)$. 

The second construction is the Alesker-Fu product \cite{ale05c}, which is a bilinear map 
\begin{displaymath}
\mathcal{V}^\infty(M) \times
\mathcal{V}^\infty(M) \to \mathcal{V}^\infty(M), (\mu_1,\mu_2) \mapsto
\mu_1 \cdot \mu_2.
\end{displaymath}
We refer to \cite{ale05c} for its construction. It is characterized by
the following properties: 
\begin{enumerate}
\item $\cdot$ is continuous and linear in both variables;
\item if $\rho:N \to M$ is a diffeomorphism and $\mu_1, \mu_2 \in
  \mathcal{V}^\infty(M)$, then 
\begin{displaymath}
\rho^*(\mu_1 \cdot
  \mu_2)=\rho^* \mu_1 \cdot \rho^*\mu_2;
\end{displaymath}
\item if $m_1,m_2$ are smooth measures on an $n$-dimensional vector
  space $V$, $A_1,A_2 \in \mathcal{K}(V)$ convex
  bodies with strictly convex smooth boundary and if $\mu_i \in
  \mathcal{V}^\infty(V),i=1,2$ is defined by  
\begin{equation}
\mu_i(K)=m_i(K+A_i), \quad K \in \mathcal{K}(V), \label{simple_vals_rn}
\end{equation}
then 
\begin{displaymath}
\mu_1 \cdot \mu_2(K)=m_1 \times m_2 ( \Delta(K)+A_1 \times A_2),
\end{displaymath}
where $\Delta:V \to V \times V$ is the diagonal embedding. 
\end{enumerate}  

Our first main theorem is the following relation between Alesker-Fu product,
integration functional, Euler-Verdier involution and Rumin
differential. 

\begin{Theorem} \label{main_theorem}
Let $\mu_1 \in \mathcal{V}^\infty(M)$ be represented by
$(\omega_1,\phi_2)$; let $\mu_2 \in \mathcal{V}_c^\infty(M)$ be
represented by $(\omega_2,\phi_2)$. Then 
\begin{equation} \label{main_equation}
\int \mu_1 \cdot \sigma \mu_2 = (-1)^n \int_{S^*M} \omega_1 \wedge
(D\omega_2+\pi^*\phi_2) + \int_M \phi_1 \wedge \pi_* \omega_2. 
\end{equation} 
\end{Theorem} 

Let us call the pairing 
\begin{align} 
\mathcal{V}^\infty(M) \times \mathcal{V}_c^\infty(M) & \to \R \nonumber\\
(\mu_1,\mu_2) & \mapsto \int \mu_1 \cdot \mu_2=:\langle
\mu_1,\mu_2\rangle \label{alesker_pairing}
\end{align}
the {\it Alesker-Poincar\'e pairing}. Note that Theorem \ref{main_theorem} is equivalent to 
\begin{equation} \label{prod_compl_forms}
\langle \mu_1,\mu_2\rangle = \int_{S^*M} \omega_1 \wedge
s^*(D\omega_2+\pi^*\phi_2)+\int_M \phi_1 \wedge \pi_* \omega_2. 
\end{equation}

From Theorem \ref{main_theorem} and from the fact that the Poincar\'e pairings on $M$ and $S^*M$ are perfect, we get the following corollary (which was first proved by Alesker). 

\begin{Corollary} (\cite{ale05d}, Thm. 6.1.1))\\ \label{perfectness_pairing}
The Alesker-Poincar\'e pairing \eqref{alesker_pairing} is a perfect pairing. 
\end{Corollary}

Some more operators on
$\mathcal{V}^\infty(M)$ were introduced in \cite{bebr07}. For this, we suppose that $M$ is a Riemannian manifold. Then $S^*M$ admits an induced metric, the {\it Sasaki metric} \cite{yaish73}.

The first operator is the derivation operator $\Lambda$ (which was denoted by $\mathfrak{L}$
in \cite{bebr07}). The metric on $S^*M$ provides a canonical
choice of $\alpha$, namely $\alpha|_{(p,[\xi])}:=\frac{1}{\|\xi\|}
\pi^*\xi$ for all $(p,[\xi]) \in S^*M$. Let $T$ be the Reeb vector field on $S^*M$
(i.e. $\alpha(T)=1$ and $\mathcal{L}_T \alpha=0$). 

If the smooth
valuation $\mu$ is represented by $(\omega,\phi)$, then $\Lambda \mu$
is by definition the valuation which is represented by
$(\mathcal{L}_T\omega+i_T \pi^*\phi,0)$. 

Let us recall the definitions of the signature operator $\mathcal{S}$
and the Laplace operator $\Delta$.
Let $*$ be the Hodge star acting on $\Omega^*(S^*M)$. Let $\mu \in \mathcal{V}^\infty(M)$ be represented by
$(\omega,\phi)$. Then $\mathcal{S}\mu$ is defined as the valuation
which is represented by $(*(D\omega+\pi^*\phi),0)$. 

The Laplace operator $\Delta$ is defined as $\Delta:=(-1)^n \mathcal{S}^2$.   

Our second main theorem shows that these operators fit well into
Alesker's theory. In fact, they are formally self-adjoint with respect
to the Alesker-Poincar\'e pairing. 

\begin{Theorem} \label{self_adjointness_theorem}
For valuations $\mu_1 \in \mathcal{V}^\infty(M)$ and $\mu_2 \in
\mathcal{V}_c^\infty(M)$, the following equations hold:
\begin{align}
\langle \Lambda \mu_1,\mu_2 \rangle & = \langle \mu_1,\Lambda
\mu_2\rangle  \label{symmetry_lambda}\\
\langle \mathcal{S}\mu_1,\mu_2\rangle & =\langle \mu_1,\mathcal{S}\mu_2\rangle \label{S_symmetry}\\
\langle \Delta \mu_1,\mu_2\rangle & = \langle \mu_1,\Delta
\mu_2\rangle. \label{laplace_symmetry}
\end{align}
\end{Theorem}

We will apply these theorems in the study of the integral geometry of
$SU(2)$. This group acts on the quaternionic line $\mathbb{H}$. In this setting, it is more natural to work with the space
$\mathcal{K}(\mathbb{H})$ of convex sets instead of manifolds with
corners. By Prop. 2.6. of \cite{ale06}, there is no
loss of generality in doing so.   

It was shown by Alesker \cite{ale04} that the space of
$SU(2)$-invariant and translation invariant valuations on the
quaternionic line $\mathbb{H}$ is of dimension $10$. For each purely
complex number $u$ of norm $1$, let $I_u$ be the complex structure
given by multiplication from the right with $u$ and $\mathbb{CP}^1_u$
the corresponding Grassmannian of complex lines (with its unique
$SU(2)$-invariant Haar measure). Alesker defined a valuation $Z_u$ by 
\begin{equation}
Z_u(K):=\int_{\mathbb{CP}^1_u} \vol(\pi_L(K)) dL, \quad K \in \mathcal{K}(\mathbb{H}). 
\end{equation}
He showed that
$Z_i,Z_j,Z_k,Z_\frac{i+j}{\sqrt{2}},Z_\frac{i+j}{\sqrt{2}},Z_\frac{i+j}{\sqrt{2}}$, together with Euler characteristic $\chi$, the volume $\vol$ and the intrinsic volumes $\vol_1,\vol_3$
form a basis of $\Val^{SU(2)}$. 
Following a suggestion of Fu, we will state the kinematic formula using a more symmetric choice. Noting that $Z_{u}=Z_{-u}$ for all $u \in S^2$, the $12$ vertices $\pm u_i,i=1,\ldots,6$ of an icosahedron on $S^2$ define $6$ valuations $Z_{u_i},i=1,\ldots,6$.

We endow $SU(2)$ with its Haar measure and the semidirect product
$\overline{SU(2)} = SU(2) \ltimes \mathbb{H}$ with the product
measure. Let $\vol_k$ denote the $k$-dimensional intrinsic volume \cite{klro97}. 

\begin{Theorem} (Principal kinematic formula for $SU(2)$) \label{main_theorem_su2}\\
Let $K,L \in \mathcal{K}(\mathbb{H})$. Then 
\begin{multline*} 
\int_{\overline{SU(2)}} \chi(K \cap \bar g L) d\bar g= \chi(K) \vol(L)+\frac{4}{3\pi} \vol_1(K) \vol_3(L)+ \\
+ \frac{17}{4} \sum_{i=1}^6 Z_{u_i}(K)Z_{u_i}(L)-\frac{3}{4}\sum_{1 \leq i \neq j \leq 6} Z_{u_i}(K)Z_{u_j}(L)+\\
\frac{4}{3\pi} \vol_3(K)\vol_1(L)+\vol(K)\chi(L).
\end{multline*}
\end{Theorem} 

This theorem implies and generalizes the Poincar\'e formulas of Tasaki \cite{tas03}, (which contained an error in some constant) as we will explain in the last section.  

\subsubsection{Acknowledgements} I wish to thank Joseph Fu for illuminating discussions, in particular for the beautiful idea of using the vertices of an icosahedron as basis for $\Val_2^{SU(2)}$. 

\section{The Alesker-Poincar\'e pairing in terms of forms}

In order to prove Theorem \ref{main_theorem} and Theorem
\ref{perfectness_pairing}, we will need three lemmas which are of independent interest. 

\begin{Lemma} (Partition of unity for valuations, \cite{ale05d},
  Prop. 6.2.1)\\ \label{lemma_partition_unity}
Let $M=\cup_i U_i$ be a locally finite open cover of $M$. Then there
exist valuations $\mu_i \in \mathcal{V}^\infty(M)$ such that $\spt \mu_i
\subset U_i$ and 
\begin{displaymath}
\sum_i \mu_i=\chi.
\end{displaymath} 
\end{Lemma}
\proof
Let $1=\sum_i f_i$ a partition of unity subordinate to $M=\cup_i
U_i$. We represent $\chi$ by $(\omega,\phi)$ and let $\mu_i$ be the
valuation represented by $(\pi^*f_i \wedge \omega,f_i \phi)$. 
\endproof

By inspecting the proof of Theorem \ref{bernig_broecker} (which uses a
local variational argument), one gets the following lemma. 
 
\begin{Lemma} \label{small_supports}
Let $\omega \in \Omega^{n-1}(S^*M)$, $\phi \in \Omega^n(M)$ and define
the smooth valuation $\mu$ by \eqref{smooth_vals}. Then 
\begin{displaymath}
\spt D\omega+\pi^*\phi \subset \pi^{-1}(\spt \mu) \text{ and } \spt \pi_*\omega
\subset \spt \mu.
\end{displaymath}
\end{Lemma}

\begin{Lemma} \label{representation_compact_support}
Let $\mu \in \mathcal{V}^\infty(M)$ be compactly supported. Then $\mu$
can be represented by a pair $(\omega,\phi) \in \Omega^{n-1}(S^*M)
\times \Omega^n(M)$ of compactly supported forms. 
\end{Lemma}
\proof
We suppose $M$ is non-compact (otherwise the statement is
trivial). Let $\mu$ be represented by a pair $(\omega',\phi')$. Then
$D\omega'+\pi^*\phi'$ is compactly supported. Since $H_c^n(S^*M)=\R$,
there exists a compactly supported form $\phi \in \Omega^n_c(M)$ such that
\begin{displaymath}
[D\omega'+\pi^*\phi']=[\pi^* \phi] \in H^n_c(S^*M).
\end{displaymath} 
In other words, there is a compactly supported form $\omega \in
H^{n-1}_c(S^*M)$ such that
$d\omega=D\omega=D\omega'+\pi^*\phi'-\pi^*\phi$. By Theorem
\ref{bernig_broecker}, the pair
$(\omega,\phi)$ represents $\mu$ up to a multiple of $\chi$. Since the
valuation represented by $(\omega,\phi)$ and the valuation $\mu$ are both compactly
supported, whereas $\chi$ is not, they have to be the same.     
\endproof

\proof[Proof of Theorem \ref{main_theorem}]
Note first that the right hand side of \eqref{main_equation} is
well-defined: since $\mu_2$ is compactly supported, the same holds
true for
$D\omega_2+\pi^*\phi_2$ and $\pi_* \omega_2$ by Lemma
\ref{small_supports}. 

Next, both sides of \eqref{main_equation} are linear in $\mu_1$ and
$\mu_2$. Using Lemma \ref{lemma_partition_unity}, we may
therefore assume that the supports of $\mu_1$ and $\mu_2$ are
contained in the support of a coordinate chart. Since the Alesker-Fu product, the Euler-Verdier involution and the integration functional are natural with respect to diffeomorphisms, it suffices to prove \eqref{main_equation} in the case where $M=V$ is a real vector space of dimension $n$. 

Let us first suppose that $\mu_1$ and $\sigma \mu_2$ are of the
type \eqref{simple_vals_rn}. We thus have $\mu_1(K)=m_1(K+A_1)$ and $\sigma \mu_2(K)=m_2(K+A_2)$
with smooth measures $m_1,m_2$ and smooth convex bodies $A_1,A_2$ with
strictly convex boundary. 

The left hand side of \eqref{main_equation} is given by 
\begin{align}
\int \mu_1 \cdot \sigma \mu_2 & =m_1 \times m_2(\Delta(V)+A_1 \times
A_2)  \nonumber \\
& = \int_V m_1((\Delta V+A_1 \times A_2) \cap V \times \{x\})
dm_2(x)  \nonumber \\
& = \int_V m_1(x-A_2+A_1) dm_2(x)  \nonumber \\
& = \int_V \mu_1(x-A_2) dm_2(x)  \nonumber \\
& = \int_V \cnc(x-A_2)(\omega_1) dm_2(x) + \int_V \left(\int_{x-A_2}
  \phi_1\right) dm_2(x). \label{product_and_forms_eq}
\end{align}

Let $A \in \mathcal{K}(V)$ be smooth with strictly convex
boundary. Its support function is defined by 
\begin{align*}
h_A:V^* & \to \R\\
\xi & \mapsto \sup_{x \in A} \xi(x).
\end{align*}
Note that $h_A$ is homogeneous of degree $1$ and that
$h_{-A}(\xi)=h_A(-\xi)$. 

Define the map $G_A:S^*V \to S^*V, (x,[\xi]) \mapsto (x+d_{\xi}h_A,[\xi])$ (since $h_A$ is
homogeneous of degree $1$, $d_{\xi}h_A \in V^{**}=V$ only depends on
$[\xi]$). $G_A$ is an orientation preserving diffeomorphism of $S^*V$. 

It is easy to show (\cite{be06}, \cite{befu06}) that for $X \in \mathcal{K}(V)$ 
\begin{equation} \label{minkowski_sum}
\cnc(X+A)=(G_A)_* \cnc(X). 
\end{equation}

We next compute that for all $(x,[\xi]) \in S^*V$ 
\begin{multline} \label{ga_and_g-a}
G_A \circ s(x,[\xi])=G_A(x,[-\xi])=(x+d_{-\xi}h_A,[-\xi])=(x-d_\xi
h_{-A},[-\xi])\\
= s(x-d_\xi h_{-A},[\xi])=s \circ G_{-A}^{-1} (x,[\xi]).
\end{multline}

Let $\kappa_2 \in \Omega^n(V)$ be the form representing the measure
$m_2$. The first term in \eqref{product_and_forms_eq} is equal to 
\begin{align}
\int_V \cnc(x-A_2)(\omega_1) dm_2(x) & = \int_V
\cnc(\{x\})(G_{-A_2}^*\omega_1) dm_2(x) \nonumber \\
& = \int_V \pi_*(G_{-A_2}^*\omega_1) \wedge \kappa_2  \nonumber \\
& = \int_{S^*V} G^*_{-A_2}\omega_1 \wedge \pi^*\kappa_2  \nonumber \\
& = \int_{S^*V} \omega_1 \wedge  (G_{-A_2}^{-1})^* \pi^* \kappa_2. \label{forms_with_G}
\end{align} 

By \eqref{minkowski_sum} we have $(-1)^n Ds^*\omega_2+ (-1)^n \pi^*\phi_2=G_A^* \pi^*
\kappa_2$. Applying $s^*$ to both sides and using \eqref{ga_and_g-a}, we get 
\begin{displaymath}
(-1)^n (D\omega_2+\pi^* \phi_2)= s^* G_{A_2}^* \pi^* \kappa_2=(G_{-A_2}^{-1})^*
\pi^* \kappa_2. 
\end{displaymath}
Hence \eqref{forms_with_G} equals $(-1)^n \int_{S^*V}
\omega_1 \wedge (D\omega_2+\pi^*\phi_2)$, which is the first term in
\eqref{main_equation}. 

By Fubini's theorem, the second term in \eqref{product_and_forms_eq} equals 
\begin{displaymath}
\int_V \left(\int_{x-A_2}
  \phi_1\right) dm_2(x) = \int_V m_2(y+A_2) \phi_1(y)= \int_V \sigma \mu_2(\{y\}) \phi_1(y).  
\end{displaymath}

For $y \in V$, we have $s_* \cnc(\{y\})=(-1)^n \cnc(\{y\})$, since the
antipodal map on $S^{n-1}$ is orientation preserving precisely if $n$
is even. Hence 
\begin{displaymath}
\sigma \mu_2(\{y\})=\pi_*\omega_2(y).
\end{displaymath}

The second term in \eqref{product_and_forms_eq} thus equals $\int_V
\phi_1 \wedge \pi_*\omega_2$, which corresponds to the second term in
\eqref{main_equation}.  

This finishes the proof in the case where $\mu_1$ and $\sigma \mu_2$ are of type \eqref{simple_vals_rn}. By linearity of both sides, \eqref{main_equation} holds true for linear combinations of such valuations. Given arbitrary $\mu_1 \in \mathcal{V}^\infty(M)$ and $\mu_2 \in \mathcal{V}^\infty_c(M)$, we find sequences $\mu_1^j \in \mathcal{V}^\infty(M)$ and $\mu_2^j \in \mathcal{V}^\infty_c(M)$ such that $\mu_1^j \to \mu_1$ and $\mu_2^j \to \mu_2$ and such that $\mu_1^j$ and $\sigma \mu_2^j$ are linear combinations of valuations of type \eqref{simple_vals_rn} (compare \cite{ale05a} and \cite{ale05b}). 

By definition of the topology on $\mathcal{V}^\infty(M)$ (see Section 3.2 of \cite{ale05b}) and the open mapping theorem, there are sequences $(\omega_1^j,\phi_1^j)$ and $(\omega_2^j,\phi_2^j)$ representing $\mu_1^j, \mu_2^j$ and converging to $(\omega_1,\phi_1),(\omega_2,\phi_2)$ in the $C^\infty$-topology. 
By what we have proved, 
\begin{displaymath}
\int \mu_1^j \cdot \sigma \mu_2^j = (-1)^n \int_{S^*M} \omega_1^j \wedge
(D\omega_2^j+\pi^*\phi_2^j) + \int_M \phi_1^j \wedge \pi_* \omega_2^j
\end{displaymath}
for all $j$. Letting $j$ tend to infinity, Equation \eqref{main_equation} follows. 
\endproof

\section{Self-adjointness of natural operators}

\proof[Proof of Theorem \ref{self_adjointness_theorem}]

Note first the following equation:
\begin{equation}
\langle \sigma \mu_1,\mu_2\rangle =(-1)^n \langle \mu_1,\sigma
\mu_2\rangle. \label{sigma_symmetry}
\end{equation}
 
This equation is immediate from \eqref{prod_compl_forms} and
the fact that $s:S^*M \to S^*M$ preserves orientation if and only if $n$
is even. 

Let $\mu_i$ be represented by $(\omega_i,\phi_i)$. By Lemma \ref{representation_compact_support} we
may suppose that $\omega_2$ and $\phi_2$ are compactly supported. 

$\Lambda \mu_i$ is represented by $\xi_i:=i_T (D\omega_i+\pi^*\phi_i)$. Since
$D\omega_i+\pi^*\phi_i=\alpha \wedge \xi_i$, we get 
\begin{align}
\langle \Lambda \mu_1,\sigma \mu_2\rangle & = (-1)^n \int_{S^*M}
\xi_1 \wedge (D\omega_2+\pi^*\phi_2) \nonumber \\
 & = (-1)^{n} \int_{S^*M} \xi_1 \wedge \alpha \wedge \xi_2 \nonumber \\
& = - \int_{S^*M} \xi_2 \wedge (D\omega_1+\pi^*\phi_1) \nonumber \\
& = (-1)^{n+1} \int_{S^*M} \omega_1 \wedge D\xi_2 - \int_M \phi_1 \wedge \pi_*
\xi_2 \nonumber \\
& = -\langle \mu_1,\sigma \Lambda \mu_2\rangle. \label{relations_for_lambda}
\end{align}

Since $D$ and $s^*$ commute and since $i_T \circ s^*=-s^* \circ i_T$, it is easily checked that 
\begin{equation} \label{commutator_lambda_sigma}
\Lambda \circ \sigma=-\sigma \circ \Lambda. 
\end{equation} 
Now \eqref{symmetry_lambda} follows from \eqref{relations_for_lambda} and \eqref{commutator_lambda_sigma}. 

Let us next prove \eqref{S_symmetry} (\eqref{laplace_symmetry} is an
immediate consequence). 

By Lemma
\ref{representation_compact_support} we may suppose that
$\omega_2$ and $\phi_2$ have compact support. Then 
\begin{align}
\langle \mu_1,\sigma \mathcal{S} \mu_2\rangle & = (-1)^n \int_{S^*M}
\omega_1 \wedge D*(D\omega_2+\pi^*\phi_2) \nonumber \\ 
& \quad + \int_M
\phi_1 \wedge \pi_**(D\omega_2+\pi^*\phi_2) \nonumber \\
& = \int_{S^*M} (D\omega_1+\pi^*\phi_1) \wedge *(D\omega_2 +
\pi^*\phi_2) \nonumber \\ 
& = \int_{S^*M} *(D\omega_1+\pi^*\phi_1) \wedge (D\omega_2 +
\pi^*\phi_2) \nonumber \\
& = \int_{S^*M} (-1)^n s^* *(D\omega_1
+ \pi^*\phi_1) \wedge s^*(D\omega_2+\pi^*\phi_2)   \nonumber \\
& = \langle \sigma \mathcal{S}\mu_1,\mu_2\rangle.
\label{sigmaS_equation} 
\end{align}

Since $s$ changes the orientation of $S^*M$ by $(-1)^n$, we get $s^* \circ *=(-1)^{n} * \circ s^*$ on $\Omega^*(S^*V)$. 
It follows that $\sigma \circ \mathcal{S}=(-1)^{n}\mathcal{S} \circ
\sigma$. Therefore \eqref{S_symmetry} follows from \eqref{sigma_symmetry} and
\eqref{sigmaS_equation}.
\endproof

Alesker defined the space $\mathcal{V}^{-\infty}(M)$ of {\it
  generalized valuations} on $M$ by 
\begin{displaymath}
\mathcal{V}^{-\infty}(M):=\left(\mathcal{V}_c^\infty(M)\right)^*,
\end{displaymath}
where the star means the topological dual. This space is endowed with
the weak topology. 
By the perfectness of the Alesker-Poincar\'e pairing, there is a
natural dense embedding $\mathcal{V}^\infty(M) \hookrightarrow
\mathcal{V}^{-\infty}(M)$.

\begin{Corollary}
Let $M$ be a Riemannian manifold. Each of the operators
$\Lambda,\mathcal{S},\Delta$ acting on $\mathcal{V}^\infty(M)$ admits
a unique continuous extension to $\mathcal{V}^{-\infty}(M)$. 
\end{Corollary}

\proof
Uniqueness of the extension is clear, since $\mathcal{V}^\infty(M)$ is
dense in $\mathcal{V}^{-\infty}(M)$. We let $\Lambda$ act on
$\mathcal{V}^{-\infty}$ by $\Lambda \xi(\mu):=\xi(\Lambda \mu)$. By
Theorem \ref{self_adjointness_theorem}, this is consistent with the
embedding of $\mathcal{V}^\infty(M)$ into $\mathcal{V}^{-\infty}(M)$
and we are done. The cases of $\mathcal{S}$ and $\Delta$ are similar.     
\endproof 

\section{The translation invariant case}

From now on, $V$ will denote an oriented $n$-dimensional real vector
space. We will consider valuations on the space
$\mathcal{K}(V)$ of compact convex
sets (i.e. convex valuations).    

A convex valuation $\mu$ on $V$ is called translation invariant, if $\mu(x+K)=\mu(K)$ for all $K \in \mathcal{K}(V)$ and all $x \in V$. 

A translation invariant convex valuation $\mu$ is said to be of degree
$k$ if $\mu(tK)=t^k \mu(K)$ for $t>0$ and $K \in \mathcal{K}(V)$. By
$\Val_k(V)$ we denote the space of translation invariant convex valuations of degree $k$. A valuation $\mu$ is even if $\mu(-K)=\mu(K)$ and odd if $\mu(-K)=-\mu(K)$, the corresponding spaces will be denoted by a superscript $+$ or $-$. 

In \cite{mcmull77} it is shown that the space of translation invariant valuations can be written as a direct sum 
\begin{displaymath} 
\Val(V)=\bigoplus_{k=0}^n \Val_k(V). 
\end{displaymath}
Each space $\Val_k(V)$ splits further as $\Val_k(V)=\Val_k^+(V) \oplus \Val_k^-(V)$. 

The spaces $\Val_0(V)$ and $\Val_n(V)$ are both $1$-dimensional (generated by $\chi$ and a Lebesgue measure respectively). For $\mu \in \Val(V)$, we denote by $\mu_n$ its component of degree $n$. 

Let us prove the following version of Theorem \ref{main_theorem} in the translation invariant case. 
\begin{Theorem} \label{main_theorem_translation_invariant}
Let $\mu_1,\mu_2 \in \Val^{sm}(V)$ be represented by translation
invariant forms $(\omega_1,\phi_1)$, $(\omega_2,\phi_2)$
respectively. Then $(\mu_1 \cdot \sigma \mu_2)_n$ is represented by the $n$-form 
\begin{displaymath} 
(-1)^n \pi_*(\omega_1 \wedge (D\omega_2+\pi^*\phi_2)) +\phi_1 \wedge \pi_*\omega_2 \in \Omega^n(V). 
\end{displaymath}
\end{Theorem}

\proof
The proof is similar to that of Theorem \ref{main_theorem}. Fix a Euclidean metric on $V$. For $R>0$, let $B_R$ denote the ball of radius $R$, centered at the origin. Let us suppose that $\mu_1(K)=\vol(K+A_1)$ and $\sigma \mu_2(K)=\vol(K+A_2)$ for all $K \in\mathcal{K}(V)$. Then 
\begin{align*}
\mu_1 \cdot \sigma \mu_2(B_R)& =\vol_{2n}(\Delta(B_R)+A_1 \times A_2)\\
& = \int_{B_R} \vol(x-A_2+A_1)dx+o(R^n)\\
& = \int_{B_R} \mu_1(x-A_2) +o(R^n)\\
& = \int_{B_R} \cnc(x-A_2)(\omega_1) dx + \int_{B_R} \int_{x-A_2} \phi_1 dx + o(R^n).
\end{align*}

The first term is given by 
\begin{align*}
\int_{B_R} \cnc(x-A_2)(\omega_1) dx & = \int_{B_R} \pi_* (G_{-A_2}^*\omega_1)dx+o(R^n)\\
& = \int_{B_R \times S^*(V)} G_{-A_2}^* \omega_1 \wedge \pi^*(dx) +o(R^n)\\
& = \int_{G_{-A_2}(B_R \times S^*(V))} \omega_1 \wedge (G_{-A_2}^{-1})^* \pi^* dx_2+o(R^n)\\
& = (-1)^n \int_{B_R \times S^*(V)} \omega_1 \wedge (D\omega_2+\pi^* \phi_2) +o(R^n)\\
& = (-1)^n \int_{B_R} \pi_*( \omega_1 \wedge (D\omega_2+\pi^* \phi_2))+o(R^n).
\end{align*}
  
The second term yields 
\begin{align*}
\int_{B_R} \int_{x-A_2} \phi_1 dx & = \int_V \vol((y+A_2) \cap B_R) \phi_1(y)\\
& = \int_{B_R} \vol(y+A_2) \phi_1(y)+o(R^n)\\
& = \int_{B_R} \mu_2(\{y\}) \phi_1(y)+o(R^n)\\
& = \int_{B_R} \phi_1 \wedge \pi_*\omega_2 +o(R^n).
\end{align*}

Therefore we obtain 
\begin{align*}
(\mu_1 \cdot \sigma \mu_2)_n & =
\lim_{R \to \infty} \frac{1}{R^n} \mu_1 \cdot \sigma \mu_2(B_R)\\
& = \lim_{R \to \infty}   \frac{1}{R^n} \int_{B_R}  (-1)^n \pi_*( \omega_1 \wedge (D\omega_2+\pi^* \phi_2))+\phi_1 \wedge \pi_*\omega_2.  
\end{align*}
This finishes the proof of Theorem \ref{main_theorem_translation_invariant} in the case where $\mu_1,\sigma \mu_2$ are of type $K \mapsto \vol(K+A)$. Using linearity of both sides, it also hold for linear combinations of such valuations. Since they are dense in $\Val(V)$ (by Alesker's solution of McMullen's conjecture \cite{ale01}), Theorem \ref{main_theorem_translation_invariant} is true in general. 
\endproof

Let us next suppose that $V$ is endowed with a Euclidean product. We can identify $\Val_n(V)$ with $\R$ by sending $\vol$ to $1$. We get a symmetric bilinear form (called Alesker pairing)  
\begin{align*}
\Val^{sm}(V) \times \Val^{sm}(V) &\to \R; \\
(\mu_1,\mu_2) & \mapsto \langle \mu_1,\mu_2\rangle:=(\mu_1 \cdot \mu_2)_n.
\end{align*} 

\begin{Corollary} \label{self_adjointness_theorem_translation_invariant}
For valuations $\mu_1,\mu_2 \in \Val^{sm}(V)$, the following equations hold:
\begin{align*}
\langle \Lambda \mu_1,\mu_2 \rangle & = \langle \mu_1,\Lambda
\mu_2\rangle\\  
\langle \mathcal{S}\mu_1,\mu_2\rangle & =\langle \mu_1,\mathcal{S}\mu_2\rangle \\ 
\langle \Delta \mu_1,\mu_2\rangle & = \langle \mu_1,\Delta
\mu_2\rangle. 
\end{align*}
\end{Corollary}

\proof
Analogous to the proof of Theorem \ref{self_adjointness_theorem}. 
\endproof

\section{Kinematic formulas and Poincar\'e formulas}

\subsection{Kinematic formulas}

In this section, we suppose that $G$ is a
subgroup of $O(V)$ acting transitively on the unit sphere. By a
result of Alesker \cite{ale04}, the space of translation invariant and
$G$-invariant valuations $\Val^G$ is a finite-dimensional vector
space.  

Let $\phi_1,\ldots,\phi_N$ a basis of $\Val^G$. Suppose we have a
kinematic formula 
\begin{displaymath}
\int_{\bar G} \chi(K \cap \bar gL)d\bar g=\sum_{i,j=1}^N
c_{i,j}\phi_i(K)\phi_j(L). 
\end{displaymath}
Here and in the following, $G$ is endowed with its Haar measure and
$\bar G:=G \ltimes V$ with the product measure. 

Set   
\begin{displaymath}
k_G(\chi):=\sum_{i,j=1}^N c_{i,j}\phi_i \otimes \phi_j \in \Val^G \otimes \Val^G=\Hom(\Val^G,\Val^{G*}).
\end{displaymath} 
The Alesker pairing induces a bijective map 
\begin{displaymath}
\PD \in \Hom(\Val^G,\Val^{G*}).  
\end{displaymath} 
Fu \cite{fu04} showed that these two maps are inverse to each other: 
\begin{equation} \label{fu_equation}
k_G(\chi)=\PD^{-1}.
\end{equation} 

For further use, we give another interpretation of
\eqref{fu_equation}. Let $G$ be as above. The scalar product on the
finite-dimensional space $\Val^G$ induces a scalar product on
$\Val^{G*}$ such that $\PD$ is an isometry. 

Given $K \in \mathcal{K}(V)$, let $\mu_K \in \Val^{G*}$ be defined by 
\begin{displaymath}
\mu_K(\mu)=\mu(K), \quad \mu \in \Val^G.
\end{displaymath}

\begin{Proposition} (Principal kinematic formula)\\
Let $G$ be a subgroup of $O(V)$ acting transitively on the unit
sphere. Then for $K,L \in \mathcal{K}(V)$
\begin{displaymath} 
\int_{\bar G} \chi(K \cap \bar gL)d\bar g=\langle
\mu_K,\mu_L\rangle.
\end{displaymath}
\end{Proposition}

\proof
Let $\phi_1,\ldots,\phi_N$ be a basis of $\Val^G$. Set $g_{ij}:=\langle \phi_i,\phi_j\rangle, i,j=1,\ldots,N$. Let us denote by
$(g^{ij})_{i,j=1\ldots,N}$ the inverse matrix. 
Then 
\begin{displaymath}
\int_{\bar G} \chi(K \cap \bar gL)d\bar g  = \sum_{i,j}
g^{ij}\phi_i(K)\phi_j(L) 
 =\sum_{i,j} g^{ij}\mu_K(\phi_i)\mu_L(\phi_j)
 = \langle \mu_K,\mu_L\rangle.
\end{displaymath}
\endproof

\subsection{Klain functions}

Let us suppose additionally that $-1 \in G$, which implies that
$\Val^G \subset \Val^+$. 
 
For $0 \leq k \leq n$, the action of $G$ on $V$ induces an action on
the Grassmannian $\Gr_k(V)$. We set $\mathcal{P}_k:=\Gr_k(V)/G$ for
the quotient space. Given $u \in \mathcal{P}_k$, the space of
$k$-planes contained in $u$ admits a unique $G$-invariant probability
measure and we define $Z_u \in \Val^G$ by 
\begin{displaymath} 
Z_u(K):=\int_{L \in u} \vol(\pi_LK)dL, \quad K \in \mathcal{K}(V).
\end{displaymath}

Recall that the Klain function of an even, translation invariant
valuation $\mu$ of degree $k$ on a Euclidean vector space $V$ is the
function $\kl_\mu:\Gr_k(V) \to \R$ such that the restriction of $\mu$
to $L \in \Gr_k(V)$ is given by $\kl_\mu(L)$ times the Lebesgue
measure. An even, translation invariant valuation is uniquely
determined by its Klain function \cite{kl00}. If $M$ is a compact
$k$-dimensional submanifold (possibly with boundaries or corners), then 
\begin{displaymath}
\mu(M)=\int_M \kl_{\mu}(T_pM)dp.
\end{displaymath}

Alesker proved the existence of a duality operator (or Fourier transform) $\mathbb{F}$ on $\Val^{+,sm}$
such that $\kl_{\mathbb{F}\mu}=\kl_{\mu} \circ \perp$ for all $\mu \in
\Val^{+,sm}$. $\mathbb{F}$ is formally self-adjoint with respect to the Alesker pairing. 

\begin{Proposition} \label{klain_function_general}
Let $u,v \in \mathcal{P}_k$ and $L \in v$. Then   
\begin{equation} \label{klain_function_eq}
\kl_{Z_u}(L)=\langle \mathbb{F} Z_u,Z_v\rangle. 
\end{equation}
\end{Proposition}

\proof
Immediate from Lemma 2.2. of \cite{befu06}. 
\endproof

\begin{Lemma} \label{finiteness}
There are finitely many elements $u_1,\ldots,u_N$ such that
$Z_{u_i},i=1,\ldots,N$ is a basis of $\Val_k^G$. 
\end{Lemma}

\proof
Let $\phi_1,\ldots,\phi_N$ be a basis of $\Val_k^G$. Let $m_i$ be the push-forward of a Crofton measure for $\phi_i$ on $\Gr_k(V)$ under the projection $\Gr_k(V) \to \mathcal{P}_k$. 

By $G$-invariance of $\phi_i$, we get 
\begin{displaymath}
\phi_i(K)=\int_{\mathcal{P}_k} \int_{L \in u} \vol(\pi_LK)dL dm_i(u).
\end{displaymath}

For sufficiently close approximations of the $m_i$ by discrete measures $\sum_{j=1}^{k_i}
c_{i,j}\delta_{u_{i,j}}$ with $u_{i,j} \in \mathcal{P}_k, c_{i,j}\in \R$, the valuations $\sum_j
c_{i,j}Z_{u_{i,j}}$ form a basis of $\Val_k^G$. Hence 
$\{Z_{u_{i,j}},i=1,\ldots,N,j=1,\ldots,k_i\}$ is a finite generating set of $\Val_k^G$, from which we can extract a finite basis. 
\endproof

\subsection{Poincar\'e formulas}

Poincar\'e formulas for $G$ are special cases of the principal
kinematic formula for $G$, when $K$ and $L$ are replaced
by smooth compact submanifolds $M_1$ and $M_2$ (possibly with boundary) of complementary dimension (note that $M_1,M_2 \in
\mathcal{P}(V)$, so there is no problem in evaluating a valuation in
$M_1$ and $M_2$). Then the right hand side of the principal kinematic
formula is the ``average number'' of intersections of $M_1$ and
$\bar g M_2$. 

\begin{Proposition} (General Poincar\'e formula)\\ \label{poincare_general}
Let $M_1,M_2$ be smooth compact
submanifolds, possibly with boundaries, of complementary dimensions
$k$ and $n-k$. 
Then 
\begin{displaymath} 
\int_{\bar G} \#(M_1 \cap \bar g M_2) d\bar g= \int_{M_1 \times M_2}
\alpha(T_pM_1,T_qM_2)dp dq
\end{displaymath}
with 
\begin{align*} 
\alpha: \mathcal{P}_k \times
\mathcal{P}_{n-k} & \to \R \\
(u,v) \mapsto \langle Z_u,Z_v\rangle. 
\end{align*}
\end{Proposition}

\proof
Let $u_1,\ldots,u_N$ be such that $Z_{u_i},i=1,\ldots,N$ is a basis of
$\Val_k^G$. Let $v_1,\ldots,v_N$ be such that $Z_{v_j},j=1,\ldots,N$
is a basis of $\Val_{n-k}^G$ (note that the dimensions of these two
spaces agree by Thm. 1.2.2 in \cite{ale03b}). Setting $g_{ij}:=\langle
Z_{u_i},Z_{v_j}\rangle$ and $(g^{ij})$ for the inverse matrix, the
principal kinematic formula implies that for all $M_1$ and $M_2$ as above 
\begin{align*}
\int_{\bar G} \#(M_1 \cap \bar g M_2) d\bar g & = \sum_{i,j}
g^{ij}Z_{u_i}(M_1)Z_{v_j}(M_2) \\
& = \int_{M_1 \times M_2} \sum_{i,j} g^{i,j}
\kl_{Z_{u_i}}(T_pM_1)\kl_{Z_{v_j}}(T_qM_2)dp dq. 
\end{align*}
This shows that 
\begin{displaymath}
\alpha(u,v)=\sum_{i,j} g^{ij} \kl_{Z_{u_i}}(u)\kl_{Z_{v_j}}(v)=\langle Z_u,Z_v\rangle;
\end{displaymath}
where the last equation follows from \eqref{klain_function_eq} and the self-adjointness of $\mathbb{F}$.  
 
\endproof


\section{Kinematic formulas for $SU(2)$} 

We apply the results of the preceding section to the special case
$G=SU(2)$ acting on the quaternionic line  
\begin{displaymath} 
\mathbb{H}=\{x_1+x_2i+x_3j+x_4k:(x_1,x_2,x_3,x_4) \in \R^4\}. 
\end{displaymath}

Since this action is transitive on the unit sphere, $\Val^{SU(2)}$ is
finite dimensional and $\Val_k^{SU(2)}$ is one-dimensional except for
$k=2$. The quotient space $\mathcal{P}_2:=\Gr_2(\mathbb{H})/SU(2)$ is
the two-dimensional projective space $\mathbb{RP}^2=S^2/\{\pm 1\}$
\cite{tas03}. Following Tasaki, we denote by
$(\omega_1(L):\omega_2(L):\omega_3(L)) \in \mathbb{RP}^2$ the class of $L \in
\Gr_2(\mathbb{H})$. 

A canonical representative in the preimage of $(a:b:c) \in
\mathbb{RP}^2$ is given by the $2$-plane spanned by $1$ and
$ai+bj+ck$. 

If $u=(a,b,c) \in S^2$, then the planes in $u$ are the complex lines
for the complex structure $I_u$ which is defined by multiplication by
$u$ from the right on $\mathbb{H}$. We will therefore write $\mathbb{CP}_u$
instead of $u$. Note that $\mathbb{CP}_u=\mathbb{CP}_{-u}$. 

The following $SU(2)$-invariant and translation invariant valuations of degree $2$ were introduced by Alesker \cite{ale04}. 
\begin{Definition}
Given $u \in \mathbb{RP}^2$, set
\begin{displaymath}
Z_u(K):=\int_{\mathbb{CP}_u} \vol(\pi_L(K)) dL, \quad K \in
\mathcal{K}(\mathbb{H}). 
\end{displaymath}
\end{Definition}

\proof[Proof of Theorem \ref{main_theorem_su2}]
From Lemma \ref{finiteness} we infer that there is a finite number of points $u_1,\ldots,u_N \in \mathbb{RP}^2$
such that $Z_{u_i},i=1,\ldots,N$ form a basis of
$\Val_2^{SU(2)}$. Alesker showed that $N=6$ and that
$Z_i,Z_j,Z_k,Z_{\frac{i+j}{\sqrt{2}}},Z_{\frac{i+k}{\sqrt{2}}},Z_{\frac{j+k}{\sqrt{2}}}$
is such a basis. 

Our aim is to compute the product $\langle Z_u,Z_v\rangle$ for $u,v
\in \mathbb{RP}^2$. We will achieve it by first expressing each $Z_u$
as a smooth valuation represented by some $3$-form $\omega_u \in
\Omega^3(S^*\mathbb{H})$ and then applying Theorem
\ref{main_theorem_translation_invariant}. 

Since the metric induces a diffeomorphism between $S^*\mathbb{H}$ and
$S\mathbb{H}$, we may as well work with the latter space. The image of
the conormal cycle of a compact convex set $K$ under this
diffeomorphism is the normal cycle $\nc(K)$. 

Let us introduce several differential forms on $S\mathbb{H}$, depending on the choice of the complex structure $I_u$. We follow the notation of \cite{pa02}. 

Let $\alpha,\beta,\gamma$ be $1$-forms on $S\mathbb{H}$ which, at a point $(x,v) \in S\mathbb{H}$, equal 
\begin{align*}
\alpha(w) & =\langle v,d\pi(w)\rangle, w \in T_{(x,v)} S\mathbb{H},\\
\beta_u(w) & =\langle v,I_u d\pi(w)\rangle, w \in T_{(x,v)}S\mathbb{H},\\
\gamma_u(w) & =\langle v,I_u d\pi_2(w)\rangle, w \in T_{(x,v)}S\mathbb{H}. 
\end{align*}

Note that $\alpha$ is the canonical $1$-form (in particular independent of $u$), whereas $\beta_u$ and $\gamma_u$ depend on $u$.  

Let $\Omega$ be the pull-back of the symplectic form on $(\mathbb{H},I_u)$ to $S\mathbb{H}$, i.e. 
\begin{displaymath}
\Omega_u(w_1,w_2):=\langle d\pi(w_1),I_u d\pi(w_2)\rangle, w_1,w_2 \in T_{(x,v)}S\mathbb{H}.
\end{displaymath}

{\bf Claim:} $Z_u$ is represented by the $3$-form
\begin{displaymath}
\omega_u:=\frac{1}{8\pi} \beta_u \wedge d\beta_u + \frac{1}{4\pi} \gamma_u \wedge \Omega_u.
\end{displaymath}

Since $\omega_u$ is $U(2)$- and translation invariant and has bidegree
$(2,1)$ (with respect to the product decomposition
$S\mathbb{H}=\mathbb{H} \times S(\mathbb{H})$), it represents some
$U(2)$-invariant, translation invariant valuation $\mu_u$ of degree
$2$. Here $U(2)$ is the unitary group for the complex structure $I_u$. 

Now the space of valuations with these properties is of dimension $2$
\cite{ale03b}. It is thus enough to show that the valuation $Z_u$ and the valuation $\mu_u$ agree on the unit ball $B$ as well as on a complex disk $D_u$. 

It is clear that $Z_u(B)=\omega_2=\pi$. It was shown by Fu (compare Equation (37) in \cite{fu04},) that $Z_u(D_u)=\frac{\pi}{2}$.

By \cite{bebr07}, the derivation of a smooth translation invariant valuation $\mu$ on a finite-dimensional Euclidean vector space is given by 
\begin{displaymath}
\Lambda \mu(K)=\left.\frac{d}{dt}\right|_{t=0} \mu(K+tB). 
\end{displaymath}

It follows that, if $\mu$ is of degree $k$, then $\Lambda \mu(B)=k \mu(B)$. 
 
It is easily checked that $\mathcal{L}_T \beta=\gamma$, $\mathcal{L}_T\gamma=0$ and $\mathcal{L}_T^2 \Omega= d\gamma$, so that 
\begin{displaymath}
\mathcal{L}_T^2 \omega_u=\frac{1}{2\pi} \gamma \wedge d\gamma. 
\end{displaymath}
Note that $\gamma \wedge d\gamma$ is twice the volume form on
$S^3$, hence $\Lambda^2 \mu_u=2\pi \chi$. It follows that 
\begin{displaymath}
\mu_u(B)=\frac12 \Lambda^2 \mu_u(B)=\pi.
\end{displaymath}

The restriction of $\beta_u$ to the normal cycle of the complex disc $D_u$ clearly vanishes. $\gamma_u$ is the length element of the fibers of $\pi:\nc(D_u) \to D_u$ (which are circles), whereas $\Omega_u$ is the (pull-back of) the volume form on $D_u$. It follows that $\omega_u(D_u)=\frac{\pi}{2}$.  The claim is proved. 

Next, the Rumin operator is easily computed as 
\begin{displaymath}
D\omega_u=d\left(\omega_u+\frac{1}{8\pi} \alpha \wedge \beta_u \wedge \gamma_u-\frac{1}{8\pi} \alpha \wedge \Omega_u \right) = \frac{1}{2\pi} \alpha \wedge \beta_u \wedge d\gamma_u. 
\end{displaymath} 

From Theorem \ref{main_theorem_translation_invariant} we infer that $\mu_u \cdot \mu_v$ is represented by the $4$-form 
\begin{displaymath}
\frac{1}{16\pi^2} \pi_*((\beta_u \wedge d\beta_u + 2
\gamma_u \wedge \Omega_u) \wedge \alpha \wedge \beta_v \wedge
d\gamma_v) \in \Omega^4(\mathbb{H}).
\end{displaymath}

If $u=(a:b:c)$ and $v=(\tilde a:\tilde b:\tilde c)$, then    
\begin{multline*}
 (\beta_u \wedge d\beta_u + 2
\gamma_u \wedge \Omega_u) \wedge \alpha \wedge \beta_v \wedge
d\gamma_v\\
=2\left((a\tilde b-\tilde a b)^2+(a\tilde c-\tilde a c)^2+(b\tilde c-\tilde b c)^2+2(a\tilde a+b\tilde b+c \tilde c)^2\right) d\vol_{S\mathbb{H}}\\
=2(1+(a\tilde a+b\tilde b+c \tilde c)^2) d\vol_{S\mathbb{H}}.
\end{multline*}

It follows that
\begin{equation} \label{product_quaternionic_vals}
\langle Z_u,Z_v\rangle= \frac14 \left(1+(a\tilde a+b\tilde b+c \tilde c)^2\right).
\end{equation}

Let $\pm u_i,i=1,\ldots,6$ be the $12$ vertices of an
icosahedron $I$ on $S^2$. They induce $6$ valuations
$Z_{u_i},i=1,\ldots,6$. Since the edge length $a$ of $I$ satisfies
$\cos a=\frac{\sqrt{5}}{5}$, \eqref{product_quaternionic_vals}
implies that  

\begin{equation} \label{product_ico}
\langle Z_{u_i},Z_{u_j}\rangle = \left\{ \begin{array}{c c} \frac{1}{2} & i=j\\ \frac{3}{10} & i \neq j \end{array} \right.
\end{equation}

Theorem \ref{main_theorem_su2} follows easily from
\eqref{product_ico} and \eqref{fu_equation}. 
\endproof

The general Poincar\'e formula (Proposition \ref{poincare_general}) implies
the following (corrected version of the) 
Poincar\'e formula on the quaternionic line.  

\begin{Corollary} (Poincar\'e formula for $SU(2)$, \cite{tas03})\\
Let $M_1,M_2 \subset \mathbb{H}$ be compact smooth $2$-dimensional submanifolds. Then 
\begin{displaymath}
\int_{\overline{SU(2)}} \#(M_1 \cap \bar gM_2) d \bar g= \frac14 \int_{M_1 \times M_2} (1+A(T_pM_1,T_qM_2)) dp dq
\end{displaymath} 
with 
\begin{displaymath}
A(T_pM_1,T_qM_2)=(\omega_1(T_pM_1)\omega_1(T_qM_2)+\omega_2(T_pM_1)\omega_2(T_qM_2)+\omega_3(T_pM_1)\omega_3(T_qM_2))^2.
\end{displaymath}
\end{Corollary}


\end{document}